\documentclass[11pt,a4paper,reqno]{amsart}

\usepackage{amsfonts, amsmath, amsthm, amssymb, amscd}

\newtheorem{theorem}{Theorem}[section]
\newtheorem{lemma}[theorem]{Lemma}
\newtheorem{proposition}[theorem]{Proposition}
\newtheorem{corollary}[theorem]{Corollary}

\theoremstyle{definition}

\theoremstyle{plain}

\begin{document}
    \title[Rational Semigroup Automata]{Rational Semigroup Automata}

\maketitle

\begin{center}
    Elaine Render and Mark Kambites \\

    \medskip

    School of Mathematics, \ University of Manchester \\
    Manchester M13 9PL, \  England \\

    \medskip

    \texttt{E.Render@maths.manchester.ac.uk} \\
    \texttt{Mark.Kambites@manchester.ac.uk} \\

    \medskip


\end{center}

\begin{abstract}
We show that for any monoid $M$, the family of languages accepted by
$M$-automata (or equivalently, generated by regular valence grammars
over $M$) is completely determined by that part of $M$ which lies outside
the maximal ideal. In particular, every such family arises as a family
of languages accepted by $N$-automata where $N$ is a simple or $0$-simple
monoid. A consequence is that every such family is either exactly the
class of regular languages, contains all the blind one-counter languages,
or is the family of languages accepted by $G$-automata where $G$ is an
non-locally-finite torsion group.

We then consider a natural extension of the usual definition which permits
the automaton to utilise more of the structure of each monoid, and
additionally allows us to define $S$-automata for $S$ an arbitrary semigroup.
In the monoid case, the resulting automata are equivalent to the \textit{valence automata with
rational target sets} which arise in the theory of regulated
rewriting systems. We study these automata in the case that the
register semigroup is completely simple or completely $0$-simple,
obtaining a complete characterisation of the classes of languages
corresponding to such semigroups, in terms of their maximal subgroups. In
the process, we obtain a number of results about rational subsets of Rees
matrix semigroups which are likely to be of independent interest.
\end{abstract}

\bigskip

\section{Introduction}
Recently there has been increasing interest in finite automata augmented
with a memory register which stores at any moment in time an element of a
given monoid $M$ (or group $G$). The register is initialised with the
identity element of the monoid, modified by right multiplication by monoid
elements, and a word is accepted if and only if some computation reading
the word returns the register to the identity with the finite state control
in an accepting state. These automata are of considerable interest
both to algebraists (who know them as \textit{blind monoid automata} or
\textit{$M$-automata}) and to
computer scientists (to whom they are often \textit{extended finite automata} or
\textit{valence automata}). On the one hand, they
provide algebraic characterizations of important language classes
such as the context free languages, the recursively
enumerable languages and the blind and partially blind counter
languages \cite{Gilman96, Kambites06, Mitrana/Stiebe97, Redko/Lisovik80, Render/Kambites07}.
On the other, they give insights
into computational problems in group and monoid theory
\cite{Elston/Ostheimer04,Gilman/Shapiro98,Kambites06}. They
are also closely related to regulated rewriting systems, and in
particular to the \textit{valence grammars} introduced by Paun
\cite{Paun80}; the languages accepted by $M$-automata are exactly the
languages generated by regular valence grammars over $M$
\cite{Fernau/Stiebe02}.

While $M$-automata with $M$ a monoid appear at first sight to
provide much more flexibility that their group counterparts, the
extent to which such an automaton can fully utilise the structure of
the register monoid is somewhat limited. Indeed, if the register
ever contains an element of a proper ideal, then no sequence of
actions of the automaton can cause it to contain the identity again;
thus, the automaton has entered a ``fail'' state from which it can
never accept a word. It follows that the automaton can make
effective use of only that part of the monoid which does not lie in
a proper ideal. This observation will be formalised below in
Section~\ref{monoid}, where we show that every $M$-automaton is
equivalent to an $N$-automaton where $N$ is a \textit{simple} or
\textit{$0$-simple} monoid.

A natural way to circumvent this is to weaken the
requirement that the identity element be the sole initial and accepting
configuration of the register, and instead permit more general sets
of initial and terminal configurations. Permitting more general
terminal sets was suggested by Gilman \cite{Gilman96}, and the idea
has also recently appeared in the study of regulated rewriting systems,
where the introduction of \textit{valence grammars with target sets}
leads naturally to a corresponding notion of \textit{valence
automata with target sets} \cite{Fernau/Stiebe01,Fernau/Stiebe02}.
An additional advantage of allowing more general initial and accepting
sets is that we remove entirely the special role played by the identity,
and hence the very need for identity. The resulting automata are thus not
constrained to have register monoids, but instead can use arbitrary
semigroups.

If we are to retain the advantages of monoid automata as an elegant
and easily manipulated way of describing important language classes,
it is clearly necessary to place some kind of restriction on the
class of subsets permitted for initial and terminal configurations.
Possible choices include the finite subsets or the finitely generated
subsemigroups, but from a computational perspective, the most natural
choice seems to be the more general \textit{rational subsets} of the
semigroup. These sets, which have been the subject of intensive study
by both mathematicians and computer scientists (see for example
\cite{Berstel78,Kambites/Silva/Steinberg07, Redko/Lisovik80, Sakarovitch87, Stephen93}), are
general enough to significantly increase the descriptive power of the
automata, while remaining sufficiently well-behaved to permit the
development of a meaningful theory.

The objective of this paper, then, is to begin the systematic study of
finite automata augmented with a semigroup register with rational initial
and accepting subsets. In Section~\ref{preliminaries}, we briefly recall
some elementary definitions from the theories of formal languages, monoid
automata and rational subsets. In Section~\ref{monoid}, we motivate the more
general constructions which follow, by exhibiting a number of limitations on
the capability of conventional $M$-automata to describe language classes.

In Section~\ref{rationalautomata} we formally define rational
semigroup automata, and obtain some foundational results about these
automata and the classes of languages which they define. In
Section~\ref{reesmatrix} we turn our attention to Rees matrix
constructions. We first study the relationship between rational
subsets and Rees matrix constructions, obtaining a number of results
which may be of independent interest. Some of these results are then
combined with a classical theorem of Rees \cite{Rees40} to yield a
complete description of the classes of languages accepted by
rational $S$-automata whenever $S$ is a completely simple or
completely $0$-simple semigroup.

\section{Preliminaries}\label{preliminaries}

Firstly, we recall some basic ideas from formal language theory. Let
$\Sigma$ be a finite alphabet. We denote by $\Sigma^*$ the set
of all words over $\Sigma$ and by $\epsilon$ the empty word. Under
the operation of concatenation and with the neutral element
$\epsilon$, $\Sigma^*$ forms the \textit{free monoid} on $\Sigma$. The
set $\Sigma^* \setminus \lbrace \epsilon \rbrace$ of non-empty words
forms a subsemigroup of $\Sigma^*$, called the \textit{free semigroup}
$\Sigma^+$ on $\Sigma$.  A \emph{finite
automaton} over $\Sigma^*$ is a finite directed graph with each edge
labelled with an element of $\Sigma$, and with a distinguished
\textit{initial vertex} and a set of distinguished \textit{terminal vertices}.
A word
$w \in \Sigma^*$ is \textit{accepted} by the automaton if there exists a
directed path connecting the initial vertex with some terminal vertex labelled
cumulatively with $w$. The set of all words accepted by the
automaton is denoted $L$ or for an automaton $A$ sometimes $L(A)$,
and is called the \emph{language} accepted by $A$. A language
accepted by a finite automaton is called \emph{rational} or
\emph{regular}.

More generally, if $S$ is a semigroup then a \textit{finite
automaton over $S$} is a finite directed graph with each edge
labelled with an element of $S$, and with a distinguished initial
vertex and a set of distinguished terminal vertices.
An element $s \in S$ is accepted by the automaton if there exists
some directed path connecting the initial vertex to some terminal vertex,
the product of whose edge labels is $s$. If $S$ is a monoid then we
admit a unique \textit{empty path} at each vertex with label the identity
element; otherwise we consider only non-empty paths.
The \emph{subset accepted} is the set
of all elements accepted; a subset of $S$ which is accepted by some
finite automaton is called a \emph{rational subset}. The rational
subsets of $S$ are exactly the homomorphic images in $S$ of regular
languages.

We shall require the following result about rational subsets of
groups, which is well-known and easy to prove.

 \begin{proposition}\label{ratinv}
 Let $G$ be a group. If $X \subseteq G$ is rational then the subset
 $X^{-1} = \{ x^{-1} \mid x \in X \}$ is also rational.
 \end{proposition}

Next we recall the usual definition of a monoid automaton. Let $M$
be a monoid with identity $1$ and let $\Sigma$ be an alphabet. An
\emph{$M$-automaton} over $\Sigma$  is a finite automaton over the
direct product $M \times \Sigma^*$. We say that it accepts a word $w
\in \Sigma^*$ if it accepts $(1, w)$, that is if there exists a path
connecting the initial vertex to some terminal vertex labelled $(1,w)$.
Intuitively, we visualize an $M$-automaton as a finite automaton
augmented with a memory register which can store an element of $M$;
the register is initialized to the identity element, is modified by
right multiplication by element of $M$, and for a word to be
accepted the element present in the memory register on completion
must be the identity element. We write $F_1(M)$ for the class of all
languages accepted by $M$-automata, or equivalently for the class of
languages generated by regular valence grammars over $M$ \cite{Fernau/Stiebe02}.

\section{Monoid Automata and Simple/$0$-Simple Semigroups} \label{monoid}

The aim of this section is to show that the extent to which an
$M$-automaton can make use of the structure of a general monoid $M$
is severely limited. In particular, we formally justify and consider
the consequences of our intuitive observation, made in the
introduction, that a monoid automaton uses only that part of the
monoid which does not lie in a proper ideal.

Recall that an \textit{ideal} $I$ of a semigroup (or monoid) $S$ is
a subset $I$ of $S$ with the property that $S^1IS^1 \subseteq I$,
where $S^1$ denotes $S$ with an identity element adjoined. To each
ideal $I$ is associated a congruence $\rho_{I}$ on $S$ such that
$(s,t) \in \rho_{I}$ if and only if either $s,t \in I$ or $s=t$. The
quotient semigroup $S / \rho_{I}$, usually just denoted $S/I$, is called a \emph{Rees
quotient}, and takes the form \[S/I = \{I\} \cup \{\{x\}\mid x \in S
\setminus I\}.\] It is usual to identify $\lbrace x \rbrace$ with
$x$ for each $x \in S \setminus I$; the element $I$, which is easily
seen to be a zero element in $S / I$, is often denoted $0$.

\begin{proposition} \label{F(M/I)}
Let $I$ be a proper ideal of a monoid $M$. Then $F_1(M) = F_1(M/I)$.
\end{proposition}
\begin{proof}
Suppose $L \in F_1(M)$, and let $A$ be an $M$-automaton accepting $L$.
First notice that any path containing an edge of the
form $(x, w)$ with $x \in I$ will itself have label with first component
in $I$; in particular, since $I$ is a proper ideal, $1 \notin I$ and
such a path cannot be an accepting path. It follows that we may remove
any such edges without changing the language accepted, and so that we
may assume without
loss of generality that $A$ has no such edges.
Now for any $x_1, \dots, x_n \in M \setminus I$, it follows from the
definition of $M/I$ that $x_1 \dots x_n = 1$ in $M$ if and only if
$\lbrace x_1 \rbrace \dots \lbrace x_n \rbrace = \lbrace 1 \rbrace$
in $M / I$. Now if we let $B$ be the $(M/I)$-automaton obtained from
$A$ by replacing edge labels of the form $(x, w)$ with $(\lbrace x
\rbrace, w)$, it follows easily that $A$ has a path from the initial
vertex to a terminal vertex labelled $(1, w)$ if and only if $B$ has
a path from the initial vertex to a terminal vertex labelled $(\lbrace 1 \rbrace, w)$.
Thus, $B$ accepts the language $L$ and $L \in F_1(M/I)$.

Conversely, if $L \in F_1(M/I)$ then $L$ is accepted by some
$(M/I)$-automaton; by an argument akin to that above, we may assume
without loss of generality that $B$ has no edges labelled by the zero
element $I$. We now obtain from $B$ a new $M$-automaton $A$ by
replacing edges labels of the form $(\lbrace x \rbrace, w)$ with
$(x,w)$, and argue as above to show that $A$ accepts exactly the
language $L$, so that $L \in F_1(M)$.
\end{proof}

A semigroup or monoid is called \emph{simple} if it does not contain
any proper ideals; similarly a semigroup with a zero element $0$ is
called \emph{$0$-simple} if its only proper ideal is $\{0\}$.

\begin{corollary}\label{0-simple}
For every monoid $M$ there is a simple or $0$-simple monoid $N$ such
that $F_1(M) = F_1(N)$.
\end{corollary}
\begin{proof}
If $M$ has no proper ideals then it is simple, so we are done.
Otherwise, let $I$ be the union of all the proper ideals of $M$.
Then $I$ is an ideal and, since the identity element $1$ does not
lie in any proper ideal, $1 \notin I$ and $I$ is a proper ideal of
$M$. Setting $N = M/I$, it is easily verified that $N$ is
$0$-simple, and by Proposition~\ref{F(M/I)} we have $F_1(M) =
F_1(M/I) = F_1(N)$ as required.
\end{proof}

Corollary~\ref{0-simple} tells us that the usual theory of
$M$-automata really only involves the very restricted classes of
simple and $0$-simple monoids.

For $S$ a semigroup, we denote by $S^0$ the semigroup with a zero
element $0$ adjoined, that is, with elements $S \cup \lbrace 0
\rbrace$ where $0$ is a new symbol not in $S$, and multiplication
given by $$st = \begin{cases} \text{the $S$-product $st$} & \text{
if } s, t \in S; \\$$ $0$ & \text{ otherwise.}
\end{cases}$$
We recall the following easy result from \cite{Render/Kambites07} which says that
adjoining a zero to a monoid $M$ makes no difference to the class of languages accepted by
$M$-automata.
\begin{proposition}\label{M0M}
Let $M$ be a monoid. Then $F_1(M^0) = F_1(M)$.
\end{proposition}

It follows from standard results in semigroup theory (see
\cite[Theorem 2.54]{Clifford/Preston61}) that a simple [$0$-simple]
monoid is either a group [respectively, a group with $0$ adjoined]
or contains an embedded copy of the \textit{bicyclic monoid} $B$. In
\cite{Render/Kambites07} we observed that the languages accepted by
$B$-automata are exactly the \textit{partially blind one-counter
languages}, while those accepted by $\mathbb{Z}$-automata are
exactly the \textit{blind one-counter languages}; both of these
classes were introduced and studied by Greibach \cite{Greibach78}.
Combining this with Propositions~\ref{F(M/I)} and \ref{M0M} we
obtain now the following.

\begin{theorem}\label{GPBCL}
Let $M$ be a monoid. Then either $F_1(M) = F_1(G)$ for some group
$G$, or $F_1(M)$ contains the partially blind one-counter languages.
\end{theorem}

Recall that a group $G$ is called \textit{torsion} if every element
has finite order, and \textit{locally finite} if every finitely
generated subgroup is finite. A straightforward consequence of
Theorem~\ref{GPBCL} is the following trichotomy.

\begin{theorem}\label{F1M}
Let $M$ be a monoid. Then $F_{1}(M)$ either
\begin{itemize}
\item[(i)] is exactly the class of regular languages; or
\item[(ii)] contains the class of blind one-counter languages; or
\item[(iii)] is equal to $F_{1}(G)$ for $G$ a torsion group which
is not locally finite.
\end{itemize}
\end{theorem}
\begin{proof}
By Theorem~\ref{GPBCL} either $F_1(M)$ contains the partially blind
one-counter languages, or $F_1(M) = F_1(G)$ for some group $G$. In
the former case, since the class of partially blind one-counter
languages contains that of blind one-counter languages, it is
immediate that (ii) holds. So suppose $F_1(M) = F_1(G)$ for some
group $G$.

If $G$ is not a torsion group then it has an element of infinite order;
this element generates a subgroup isomorphic to
$\mathbb{Z}$, from which it follows that $F_1(G)$ contains the class
$F_1(\mathbb{Z})$ of blind one-counter languages and (ii) again holds. Now by
\cite[Proposition~1]{Kambites062}, every language in $F_1(G)$ is in
$F_1(H)$ for some finitely generated subgroup of $H$ of $G$. If $G$
is locally finite, then such an $H$ must be finite, and so every
language in $F_1(G)$ is regular. Conversely, $F_1(G)$ certainly contains
the regular languages, so (i) holds. There remains
only the case in which $G$ is a torsion group which is not
locally finite, in which case (iii) holds.
\end{proof}
Elder and Mintz \cite{Elder/Mintz/private} have observed that for $G$ a finitely generated
infinite torsion group, $F_1(G)$ cannot contain the blind one-counter
languages, but does always contain non-regular languages (for example, the
word problem of $G$). It follows that the three possible conditions in
Theorem \ref{F1M} are in fact mutually exclusive. The theorem is of particular
interest because torsion groups which are not locally finite are rather
rare and difficult to construct. It would be interesting to study the
language classes $F_1(G)$ corresponding to particular known examples of such groups
\cite{Golod68,Grigorchuk80,Olshanskii82}.

\section{Rational Semigroup Automata}\label{rationalautomata}

In Section~\ref{monoid}, we saw that the extent to which traditional
monoid automata can utilise the differences in structure between
groups and monoids was limited. In this section, we consider a
generalisation which allows us to utilise more of the structure of
arbitrary monoids, and indeed semigroups.

Let $S$ be a semigroup and $\Sigma$ a finite alphabet. We define a
\emph{rational $S$-automaton} over $\Sigma$ to be a finite automaton
over the direct product $S \times \Sigma^*$ together with two
rational subsets $X_0, X_1 \subseteq S$ called the \textit{initial
set} and \textit{terminal set} respectively. The automaton accepts a
word $w \in \Sigma^*$ if there exists $x_0 \in X_0$ and $x \in S$
such that $x_0 x \in X_1$, and $(x,w)$ labels a path from the
initial vertex to a terminal vertex in the automaton. For $S$ a
semigroup, we let $F_{Rat}(S)$ be the class of languages accepted by
rational $S$-automata.

Intuitively, a rational $S$-automaton is a non-deterministic finite
automaton augmented with a register which stores an element of $S$. The
register is (non-deterministically) initialised with an element of $X_0$,
and a word is accepted if there is a computation which reads it and
leaves the finite state control in an accept state and the register
containing a value from $X_1$.

Note that if $S = M$ is a monoid then a rational $M$-automaton with
initial set $\lbrace 1 \rbrace$ is just an \textit{valence
automaton over $M$ with rational target set} of the kind studied by Fernau
and Stiebe \cite{Fernau/Stiebe01} and the present authors
\cite{Render/Kambites07}. Indeed, the following proposition says that,
for $M$ a monoid, the initial set can be taken to be $\lbrace 1
\rbrace$ without loss of generality.

\begin{proposition}\label{GX}
Let $M$ be a monoid with identity $1$, and $L \subseteq \Sigma^*$ a
language. If $L \in F_{Rat}(M)$ then $L$ is accepted by a rational
$M$-automaton with initial set $\lbrace 1 \rbrace$, that is, by an
valence automaton over $M$ with rational target set.
\end{proposition}

\begin{proof}
Let $A$ be a rational monoid automaton with initial set $X_{0}
\subseteq M$ and terminal set $X_{1} \subseteq M$ which accepts the
language $L$. Let $L' \subseteq M \times \Sigma^*$ be the full
subset accepted by $A$ interpreted as an automaton over $M \times
\Sigma^*$. Since $X_0 \subseteq M$ is rational, the set
$$X_0' = \lbrace (x, \epsilon) \mid x \in X_0 \rbrace \subseteq M \times \Sigma^*$$
is rational. Now let
$$K = X_0' L' = \lbrace (x_0 x, w) \mid x_0 \in X_{0}, (x, w) \in L'\rbrace.$$
Then $w \in L$ if and only if there exists $x_0 \in X_0$ and $x \in
X$ such that $(x, w) \in L'$ and $x_0 x \in X_1$. But this is true
exactly if $(x', w) \in K$ for some $x' \in X_1$.

Now $K$ is a product of two rational subsets, and hence is a
rational subset. Let $B$ be a finite automaton over $M \times
\Sigma^*$ recognizing $K$. If we interpret $B$ as a rational
$M$-automaton with initial set $\lbrace 1 \rbrace$ and terminal set
$X_1$, it is immediate that $B$ accepts exactly the language $L$.
\end{proof}

Combining Corollary~\ref{GX} with a result of Fernau and Stiebe
\cite{Fernau/Stiebe01} we obtain the following.
\begin{theorem}\label{GXX}
If $G$ is a group then $F_{Rat}(G) = F_1(G)$.
\end{theorem}
\begin{proof}
If $L \in F_{Rat}(G)$ then by Corollary~\ref{GX}, $L$ is accepted by
a rational $G$-automaton with initial set $\lbrace 1 \rbrace$ and
some rational terminal set $X_1$, that is, by a valence automaton with
rational target set. But now by \cite[Theorem~3.5]{Fernau/Stiebe01},
$L$ is accepted by a finite valence automaton, that is, a $G$-automaton,
so that $L$ is in $F_1(G)$ as required. The converse is immediate.
\end{proof}

We now turn our attention to the relationship between rational relations
and rational semigroup automata. Let $\Omega$ and $\Sigma$ be finite alphabets,
and consider a finite automaton over the direct
product $\Omega^+ \times \Sigma^*$; it recognizes a \emph{rational
relation} $R \subseteq \Omega^+ \times \Sigma^*$. The image of a
language $L \subseteq \Omega^+$ under the relation $R$ is the set of words
$y \in \Sigma^*$ such that $(x,y) \in R$ for some $x \in L$.

If $X_0, X_1 \subseteq S$ then their \textit{difference} is the set
$$X_0^{-1} X_1 \ = \ \lbrace x \in S \mid x_0 x = x_1 \text{ for some } x_0 \in X_0, x_1 \in X_1 \rbrace.$$
We say that a subset $X \subseteq S$ is a \textit{rational set difference} if
there exist rational subsets $X_0, X_1 \subseteq S$ such that $X = X_0^{-1} X_1$.
Note that in a group, the rational set differences are exactly the
rational subsets, but in a general semigroup this does not hold. The following
statement is a semigroup analogue of \cite[Proposition 3.1]{Render/Kambites07},
which in turn generalises a well-known observation concerning $M$-automata
(see for example \cite[Proposition~2]{Kambites062}).
\begin{proposition} \label{transduction}
Let $X_0$ and $X_1$ be subsets of a semigroup $S$, and let
$L \subseteq \Sigma^*$ be a language. Then the following are equivalent:
\begin{itemize}
\item[(i)] $L$ is accepted by a $S$-automaton with initial set $X_0$
and terminal set $X_1$;
\item[(ii)] there exists a finite alphabet $\Omega$, a morphism
$\omega : \Omega^+ \to S$ and a rational relation $\rho \subseteq \Omega^+ \times \Sigma^*$ such that 
$$L = (X_0^{-1} X_1) \omega^{-1} \rho.$$
\end{itemize}
If $S$ is finitely generated then the following condition is also
equivalent to those above.
\begin{itemize}
\item[(iii)] for every finite choice of generators $\omega : \Omega^+ \to S$
for $S$, there exists a rational relation $\rho \subseteq \Omega^+ \times \Sigma^*$ such that 
$$L = (X_0^{-1} X_1) \omega^{-1} \rho.$$
\end{itemize}
\end{proposition}
\begin{proof}
To show that (i) implies (ii), suppose that $L$ is accepted by an
$S$-automaton $A$ with initial set $X_0$ and terminal set $X_1$. Choose a
finite alphabet $\Omega$ and a map $\omega: \Omega^+ \to S$ such that the image $\Omega^+ \omega$ contains every
element of $S$ which forms the first component of an edge label in
the automaton. We now obtain from $A$ a finite automaton $B$ over
$\Omega^+ \times \Sigma^*$ by replacing each edge label $(s,x)$
with $(w,x)$ for some $w \in \Omega^+$ is such that $w\omega =
s$. It is a routine exercise to verify that
$L$ is the image of $(X_0^{-1} X_1) \omega^{-1}$ under the relation
accepted by $B$.

Conversely, suppose we are given a map $\omega : \Omega^+ \to S$ and
an automaton $B$ over $\Omega^+ \times \Sigma^*$ such that $L$ is the image
under the relation accepted by $B$ of the language $(X_0^{-1} X_1) \omega^{-1}$.
We construct from $B$ a new automaton $A$
over $S \times \Sigma^*$ by applying the map $\omega$ to the first component
of each edge label. Now interpreting $A$ as an $S$-automaton with
initial set $X_{0}$ and terminal rational set $X_{1}$, it is easily
seen that $A$ accepts the lamguage $L$.

Suppose now that $S$ is finitely generated. Clearly (iii) implies
(ii). Conversely, if (i) holds then we can extend $\omega$ arbitrarily
to a finite choice of generators $\omega' : (\Omega')^+ \to S$ for
$M$, and check that we still have the desired property, so that
(iii) holds.
\end{proof}

As a corollary, we immediately obtain the following characterisation
for language classes of the form $F_{Rat}(S)$.

\begin{proposition}\label{semitrans}
Let $S$ be a semigroup and $L \subset \Sigma^*$ a language. Then the
following are equivalent.
\begin{itemize}
\item[(i)] $L \in F_{Rat}(S)$;
\item[(ii)] there exists an alphabet $\Omega$, a morphism
$\omega : \Omega^+ \to S$, a rational set difference $X \subseteq S$
and a rational relation $\rho \subseteq \Omega^+ \times \Sigma^*$
such that $L = X \omega^{-1} \rho$.
\end{itemize}
If $S$ is finitely generated then the following condition is also
equivalent to those above.
\begin{itemize}
\item[(iii)] there exists a rational set difference $X \subseteq S$ such that
for every finite choice of generators $\omega : \Omega^+ \to S$ for
$S$, there exists a rational relation $\rho \subseteq \Omega^+
\times \Sigma^*$ such that $L = X \omega^{-1} \rho$.
\end{itemize}
\end{proposition}

Note that, unlike in the monoid case \cite{Render/Kambites07}, we cannot
conclude that $F_{Rat}(S)$ is a rational cone. This is because the
composition of a rational relation in $\Omega^+ \times \Sigma^*$ with a
rational transduction from $\Sigma^*$ to another free monoid $\Gamma^*$
need not be a rational relation in $\Omega^+ \times \Gamma^*$ (although it
will be rational in $\Omega^* \times \Gamma^*$).

\section{Rees Matrix Constructions, Completely Simple and Completely $0$-Simple Semigroups}\label{reesmatrix}

In this section we apply the results of the previous sections to
obtain a description of language classes $F_{Rat}(S)$ for semigroups
$S$ belonging to the important classes of completely simple and completely
$0$-simple semigroups.

Recall that an idempotent $e$ in a semigroup is called
\textit{primitive} if for every non-zero idempotent $f$ such that
$ef = fe = f$ we have $e = f$. A semigroup is \emph{completely
simple} [respectively, completely 0-simple] if it is simple
[0-simple] and has a primitive idempotent. For more information
about completely simple and completely $0$-simple semigroups, see
\cite{Howie95}.

Now let $T$ be a semigroup, $0$ be a new symbol not in $T$ and let
$I,J$ be non-empty sets. Let $P = (P_{ji})$ be a $J \times I$ matrix
with entries in $T \cup \{ 0 \}$. We define a new semigroup with set
of elements
$$(I \times T \times J) \ \cup \ \{ 0 \}$$
and multiplication defined by
\[ (i,t,j)(i',t',j') = \left\{ \begin{array}{ll} (i,tP_{ji'}t',j') & \textrm{if} P_{ji'} \neq 0 \\
0 & \textrm{otherwise}, \end{array} \right. \] and
\[ (i,t,j)0 = 0(i,t,j) = 00 = 0.\]
It is simple to verify that this binary operation is associative; we
call the semigroup constructed in this way a \emph{Rees matrix
semigroup with zero over $T$}, and denote it $M^{0}(T;I,J;P)$. The
semigroup $T$ is called the \emph{base semigroup} and the matrix $P$
the \emph{sandwich matrix} of the construction. If $P$ contains no
zero entries then $I \times T \times J$ forms a subsemigroup of
$M^0(T; I, J; P)$, called a \emph{Rees matrix semigroup (without
zero) over $T$} and denoted $M(T; I, J; P)$.

Rees matrix semigroups play a crucial role in much of the structural
theory of semigroups. Of particular importance is the case that the
base semigroup $T$ is a group $G$. A Rees matrix semigroup with zero over a group is called
\textit{regular} if every row and every column of the sandwich
matrix contains a non-zero entry. The importance of this
construction can be seen from the following seminal result of Rees
\cite{Rees40}.

\begin{theorem}[The Rees Theorem]
Let $S = M(G; I, J; P)$ [respectively, $S=M^{0}(G;I,J;P)$] be a
[regular] Rees matrix semigroup over a group. Then $S$ is a
completely simple [respectively, completely 0-simple] semigroup.
Conversely, every completely simple [completely 0-simple] semigroup
is isomorphic to one constructed in this way.
\end{theorem}

We shall need the following proposition.

\begin{proposition} \label{subset}
Let $S = M(T;I,J;P)$ or $S = M^0(T; I, J; P)$ be a Rees matrix
semigroup with or without zero over a semigroup $T$. Let $X
\subseteq S$ be a rational subset and suppose $i \in I$ and $j \in
J$. Then the set
$$X_{ij} \ = \ \{ g \in T \mid (i,g,j) \in X \} \ \subseteq \ T$$
is a rational subset of $T$.
\end{proposition}
\begin{proof}
Let $A$ be a finite automaton over $S$ accepting the rational subset $X$,
with vertex set $Q$. Let $J'$ be the set of all $j \in J$ such that
$A$ has an edge label with third component $j$; note that $J'$ is
necessarily finite. We construct from $A$ a new finite automaton $B$
over $T$ with
\begin{itemize}
\item vertex set $(Q \times J') \cup \lbrace q_0' \rbrace$ where $q_0'$ is a new symbol;
\item start vertex $q_0'$;
\item terminal vertices $(q, j)$ such that $q$ is a terminal vertex of $A$;
\item an edge from $q_0'$ to $(q_1, j_1)$ labelled $t_1$ whenever $A$ has
an edge from the initial vertex to $q_1$ labelled $(i, t_1, j_1)$;
\item for every $j_1 \in J'$, an edge from $(q_1, j_1)$ to $(q_2, j_2$)
labelled $P_{j_1 i_2} t_2$ whenever $A$ has an edge from $q_1$ to
$q_2$ labelled $(i_2, t_2, j_2)$ with $P_{j_1 i_2} \neq 0$.
\end{itemize}

Since $J'$ is finite and $A$ has finitely many vertices and edges, we
deduce that $B$ has finitely many vertices and edges. Now we show that
the subset accepted by $B$ is exactly $X_{ij}$. Let $t \in X_{ij}$.
Then $(i,t,j) \in X$ labels a path in $A$ from the initial vertex
to some terminal vertex. Clearly this path cannot contain
edges labelled $0$, so it must have the form
\[
\mathbf{p_{0}}\xrightarrow{(i_1,t_{1},j_{1})}\mathbf{p_{1}}\xrightarrow{(i_{2},t_{2},j_{2})}\mathbf{p_{2}}\xrightarrow{(i_{3},t_{3},j_{3})}
\ldots
\xrightarrow{(i_{m-1},t_{m-1},j_{m-1})}\mathbf{p_{m-1}}\xrightarrow{(i_{m},t_{m},j_m)}\mathbf{p_{m}}
\] where $p_0$ is the initial vertex of $A$ and $p_m$ is a terminal vertex.
Since the path is labelled $(i, t, j)$ we must have
\begin{eqnarray*} (i,t,j) &=&
(i_{1},t_{1},j_{1})(i_{2},t_{2},j_{2}) \ldots (i_{m},t_{m},j_{m})
\end{eqnarray*}
so that $i_1 = i $, $j_m = j$. Now it follows easily from the construction
of $B$ that it has a path
\[
\mathbf{q_0'}\xrightarrow{t_{1}}\mathbf{(p_{1},j_{1})}\xrightarrow{P_{j_{1}i_{2}}t_{2}}\mathbf{(p_{2},j_{2})}\ldots
\mathbf{(p_{m-1},j_{m-1})}\xrightarrow{P_{j_{m-1}i_{m}}t_{m}}\mathbf{(p_{m},j)},
\]
where $(p_m, j)$ is a terminal vertex of $B$, so that $B$ accepts
$$t = t_{1} P_{j_{1}i_{2}} t_{2} P_{j_{2} i_3} \ldots P_{j_{m-1} i_m} t_m.$$
Thus $X_{ij} \subseteq L(B)$.

Conversely, assume that $t \in T$ is accepted by $B$. Then there
exists a path through $B$ from the initial vertex to some terminal
vertex labelled with $t$. It follows from the definition
of $B$ that this path must have the form
\[
\mathbf{q_0'}\xrightarrow{t_{1}}\mathbf{(p_{1},j_{1})}\xrightarrow{P_{j_{1}i_{2}}t_{2}}\mathbf{(p_{2},j_{2})}\ldots
\mathbf{(p_{m-1},j_{m-1})}\xrightarrow{P_{j_{m-1}i_{m}}t_{m}}\mathbf{(p_{m},j)},
\] where $p_m$ is a terminal vertex in $A$,
\begin{eqnarray*}
t &=& t_{1}P_{j_{1}i_{2}}t_{2}P_{j_{2}i_{3}}t_{3}\ldots
P_{j_{m-1}i_{m}}t_{m}
\end{eqnarray*}
and $A$ has a path
\[
\mathbf{p_{0}}\xrightarrow{(i,t_{1},j_{1})}\mathbf{p_{1}}\xrightarrow{(i_{2},t_{2},j_{2})}\mathbf{p_{2}}\xrightarrow{(i_{3},t_{3},j_{3})}
\ldots
\xrightarrow{(i_{m-1},t_{m-1},j_{m-1})}\mathbf{p_{m-1}}\xrightarrow{(i_{m},t_{m},j)}\mathbf{p_{m}}
\] where $p_0$ is the initial vertex of $A$.
Hence, $A$ accepts the element
\begin{eqnarray*}
(i,t_{1},j_{1})(i_{2},t_{2},j_{2})\ldots (i_{m},t_{m},j) &=& (i,
t_{1}P_{j_{1}i_{2}}t_{2}P_{j_{2}i_{3}}t_{3}\ldots
P_{j_{m-1}i_{m}}t_{m},j) \\
&=& (i,t,j).
\end{eqnarray*}
So $(i,t,j) \in X$ and hence $t \in X_{ij}$.

So  the automaton $B$ accepts exactly the set $X_{ij}$, and hence
$X_{ij}$ is a rational subset of $T$.
\end{proof}

As a corollary, we obtain a result about the intersections of
rational subsets with maximal subgroups in completely simple semigroups.
\begin{corollary}\label{ratmaxsubgroup}
Let $H$ be a maximal subgroup of a completely simple or completely
$0$-simple semigroup $S$. Let $X$ be a rational subset of $S$. Then
$X \cap H$ is a rational subset of $H$.
\end{corollary}
\begin{proof}
By the Rees theorem, we may assume that $S$ is a Rees matrix
semigroup without zero $M(G;I,J;P)$ or a regular Rees matrix
semigroup with zero $S = M^0(G; I, J; P)$ over a group $G$. It
follows easily from the definition of the Rees matrix construction
that either $H = \lbrace 0 \rbrace$ or
$$H \ = \ \lbrace (i, g, j) \mid g \in G \rbrace$$
for some $i \in I$ and $j \in J$ with $P_{ji} \neq 0$. In the former
case the result is trivial, so we assume the latter.
By Proposition~\ref{subset}, the set
$$X_{ij} \ = \ \lbrace g \in G \mid (i, g, j) \in X \rbrace \ = \ \lbrace g \in G \mid (i, g, j) \in H \cap X \rbrace$$
is a rational subset of $G$. It follows that
$$P_{ji} X_{ij} \ = \ \lbrace P_{ji} g \mid g \in X_{ij} \rbrace \ = \ \lbrace P_{ji} g \mid (i, g, j) \in X \rbrace$$
is also a rational subset of $G$. Now define a map
$$\phi \ : \ G \to H, \ g \mapsto (i, P_{ji}^{-1}g, j)$$
where $P_{ji}^{-1}$ is the inverse of $P_{ji}$ in the group $G$.
It is readily verified that $\phi$ is an isomorphism from $G$ to
$H$, and so the image
\begin{align*}
(P_{ji} X_{ij}) \phi \ &= \ \lbrace (i, P_{ji}^{-1} g, j) \mid g \in P_{ji} X_{ij} \rbrace \\
&= \ \lbrace (i, P_{ji}^{-1} P_{ji} g, j) \mid (i, g, j) \in X \rbrace \\
&= \ \lbrace (i, g, j) \mid (i, g, j) \in X \rbrace \\
&= \ X \cap H
\end{align*}
is a rational subset of $G$, as required.
\end{proof}

In a completely simple semigroup, where every element lies in a maximal
subgroup, Corollary~\ref{ratmaxsubgroup} easily yields the following
complete characterisation of rational subsets.

\begin{theorem}
The rational subsets of a completely simple semigroup are exactly
the finite unions of rational subsets of maximal subgroups.
\end{theorem}
\begin{proof}
Let $S$ be a completely simple semigroup. If $X_1, \dots, X_n$ are
rational subsets of maximal subgroups of $S$ then certainly they are
rational subsets of $S$, and so is their union. Conversely, suppose
$X$ is a rational subset of $S$. It follows easily from the Rees
theorem that $X$ lies inside a finitely generated completely simple
subsemigroup $S'$ of $S$. Now $S'$ is the union of finitely many
maximal subgroups, so $X$ is the union of its intersections with
these subgroups. By Corollary~\ref{ratmaxsubgroup} these
intersections are rational, so $X$ is a finite union of rational
subsets of maximal subgroups of $S'$. But maximal subgroups of $S'$ are
subgroups of $S$, and hence lie in maximal subgroups of $S'$. It follows that $X$
is a finite union of rational subsets of maximal subgroups of $S$,
as required.
\end{proof}

\begin{proposition}\label{ratcs}
Let $S = M(T;I,J;P)$ or $S = M^0(T; I, J; P)$ be a Rees matrix
semigroup with or without zero over a semigroup $T$, and let $P' \subseteq T$ be
the set of non-zero entries of the sandwich matrix $P$. Suppose $T = P' T$ or $T = T P'$.
Then for any $i \in I$, $j \in J$ and rational subset $X$ of $T$, the set
$$\lbrace (i, t, j) \mid t \in X \rbrace$$
is a rational subset of $S$.
\end{proposition}
\begin{proof}
By symmetry of assumption, it suffices to consider the case in which
$T = P' T$.
Let $A$ be a finite automaton over $T$ accepting $X$, with vertex
set $Q$. Let $Y \subseteq T$ be the set of edge labels in $A$, and
for every $t \in Y$, let $j_t \in J$, $i_t \in I$ and $s_t \in T$ be
such that $t = P_{j_t i_t} s_t$. Let $J' = \{ j_t \mid t \in Y \}
\cup \{ j \}$. Then $J'$ is a finite subset of $J$. We define a new
automaton $B$ over $S$ with
\begin{itemize}
\item vertex set $(Q \times J') \cup \lbrace q_0 \rbrace$ where $q_0$ is a new
symbol;
\item initial vertex $q_0$;
\item terminal vertices $(q, j)$ such that $q$ is a terminal vertex of $A$;
\item for every edge in $A$ from the start vertex to a vertex $q$ labelled
$t$, and every $j' \in J'$, an edge from $q_0$ to $(q, j')$ labelled $(i, t,
j')$;
\item for every edge in $A$ from a vertex $p$ to a vertex $q$ labelled $t$, and every $j' \in J'$,
an edge from $(p, j_t)$ to $(q, j')$ labelled $(i_t, s_t, j')$;
\end{itemize}
A routine argument, akin to that in the proof of
Proposition~\ref{subset}, shows that $B$ accepts the required subset
of $S$.
\end{proof}

Note in particular that the conditions on the sandwich matrix in the
hypothesis of Proposition~\ref{ratcs} are satisfied in the case of a
regular Rees matrix construction over a group.

Recall that the \textit{rational subset problem} for a finitely semigroup
$S$ is the algorithmic problem of deciding, given a rational subset (described
as a finite automaton over a fixed generating set for $S$) and an element
(described as a word over the same generating set), deciding whether the latter
belongs to the former. While the phrasing of the problem is dependent on the
precise choice of finite generating set, the decidability or undecidability
of the problem is independent of this choice \cite[Corollary~3.4]{Kambites/Silva/Steinberg07}, so one can
meaningfully say that the abstract semigroup $S$ has decidable or undecidable
rational subset problem.

\begin{corollary}
Let $S = M(T;I,J;P)$ or $S = M^0(T; I, J; P)$ be a finitely
generated Rees matrix semigroup with or without zero over a
semigroup $T$. If $T$ has decidable rational subset problem then $S$
has decidable rational subset problem.
\end{corollary}
\begin{proof}
We prove the statement for Rees matrix semigroups with zero. The
result for Rees matrix constructions without zero can be obtained as
an easy consequence, or proved directly using a similar method.

Let $\omega : \Omega^* \to T$ and $\sigma : \Sigma^* \to S$ be
finite choices of generators for $T$ and $S$ respectively. For every
$x \in \Sigma$ such that $x \sigma \neq 0$, suppose $x \sigma =
(i_x, g_x, j_x)$ and let $w_x \in \Omega^*$ be a word with $w_x
\omega = g_x$. For $j \in J$ and $i \in I$ such that $P_{ji} \neq 0$
let $w_{ji} \in \Omega^*$ be a word with $w_{ji} \omega = P_{ji}$.

Now suppose we are given a word $w = w_1 \dots w_n \in \Sigma^*$,
where each $w_i \in \Sigma$, and a rational subset $X$ of $S$.
Clearly, we can test whether $w$ represents $0$ and, in the case
that it does, whether $0 \in X$. Assume now that $w$ does not
represent $0$. Then
$$w \omega = (w_1 \omega) \dots (w_n \omega) = (i_{w_1}, g_{w_1} P_{j_{w_1} i_{w_2}} g_{w_2} \dots g_{w_n}, j_{w_n}).$$
Let $Y = \lbrace t \in T \mid (i_{w_1}, t, j_{w_n}) \in X \rbrace$, so that
$w \omega \in X$ if and only if
\begin{equation}\label{yeqn}
(w_{g_{w_1}} w_{j_{w_1} i_{w_2}} w_{g_{w_2}} \dots w_{g_{w_n}}) \sigma = g_{w_1} P_{j_{w_1} i_{w_2}} g_{w_2} \dots g_{w_n} \in Y.
\end{equation}
Now by Proposition \ref{subset}, $Y$ is rational and it follows moreover from
the proof that we can effectively compute an automaton for $Y$. By assumption,
we can solve the rational subset problem for $Y$, so we can decide whether
\eqref{yeqn} holds, as required.
\end{proof}

We now turn our attention to languages accepted by rational $S$-automata,
where $S$ is a Rees matrix semigroup. We begin with a lemma which simplifies
the case of Rees matrix semigroups with zero, by allowing us to restrict
attention to automata for which neither the initial set nor the terminal
set contains the zero element.
\begin{lemma}\label{nozero}
Let $S = M^0(T; I,J; P)$ be a finitely generated Rees matrix
semigroup with zero over a semigroup $T$. If $L$ is accepted by a
rational $S$-automaton, then $L$ is accepted by a rational
$S$-automaton for which neither the initial set nor the terminal set
contain $0$.
\end{lemma}
\begin{proof}
Suppose $L$ is accepted by a rational $S$-automaton $A$ with initial set
$X_0$ and terminal set $X_1$. Suppose first that $0 \in X_0$. If also
$0 \in X_1$ then we have $0 x \in X_1$ for all $x \in S$, so the
language accepted is just the set of all words $w$ such that $(x,w)$ labels
a path from the initial vertex to a terminal vertex of $A$ for
some $x \in S$. It follows that $L$ is regular, and hence lies in
$F_1(S)$. On the other hand, if $0 \notin X_1$ then there is no $x \in S$
such that $0x \in X_1$; hence we may replace the initial set $X_0$ with
$X_0 \setminus \lbrace 0 \rbrace$ without changing the language accepted.
Thus, we may assume that $0 \notin X_0$.

Clearly we can write $L = L_0 \cup L_1$ where $L_1$ is accepted by a
rational $S$-automaton with $0$ not in the initial or terminal sets,
and $L_0$ is accepted by a rational $S$-automaton with terminal set $0$. We
claim that $L_0$ is regular; it will follow that $L$ is the union of
$L_1$ with a regular language, and so can clearly be accepted by a
rational $S$-automaton without $0$ in the terminal set.

Let $\omega : \Omega^* \to S$ be a finite choice of generators for
$S$. For each $x \in \Omega$ such that $x \omega \neq 0$ suppose $x
\omega = (i_x, g_x, j_x)$. Now let $K$ be the set of all words
representing elements of the initial set of $A$, and let $K'$ be the
(necessarily finite) set of all final letters of words in $K$. It is
easily seen that the language
$$\lbrace v \in \Omega^* \mid (wv) \omega = 0 \text{ for some } w \in K \rbrace$$
is regular. Indeed, it consists of all words which
\begin{itemize}
\item contain a generator representing zero; or
\item contain consecutive generators $x$ and $y$ with $P_{j_x i_y} = 0$; or
\item start with a generator $y$ with $P_{j_x i_y} = 0$ for some $x \in K'$
\end{itemize}
and so can be easily described by a regular expression.
It now follows from Proposition~\ref{transduction} that $L_0$ is
a rational transduction of the above regular language and hence is itself
regular. This completes the proof.
\end{proof}

We are now ready to prove the main theorem of this section, the essence of
which is that rational $S$-automata where $S$ is a completely simple or
completely $0$-simple semigroup are no more powerful than $G$-automata
where $G$ is the maximal subgroup of $S$.

\begin{theorem}
Let $S$ be a completely simple or completely 0-simple semigroup with
maximal non-zero subgroup $G$. Then
$$F_{Rat}(S) = F_{Rat}(G) = F_1(G).$$
\end{theorem}
\begin{proof}
That $F_{Rat}(G) = F_1(G)$ is Theorem~\ref{GXX}, while
the inclusion $F_{Rat}(G) \subseteq F_{Rat}(S)$ is immediate. Hence,
we need only prove that $F_{Rat}(S) \subseteq F_{Rat}(G)$.
 It follows easily from the Rees
theorem that every completely simple semigroup $S$ embeds in a
completely 0-simple semigroup $S'$ with the same maximal
non-zero subgroup, so that $F_{Rat}(S) \subseteq F_{Rat}(S')$. Hence, it
suffices to prove the result in the case that $S$ is completely
$0$-simple.

Suppose, then, that $S$ is completely $0$-simple. By the Rees
theorem, we may assume that $S$ is a regular Rees matrix semigroup
$M^0(G^0;I,J;P)$ where $G$ is a group. Suppose now that a language
$L \subseteq \Sigma^*$ lies in $F_{Rat}(S)$. Let $A$ be a rational
$S$-automaton accepting $L$, with initial rational set $X_{0}
\subseteq S$ and terminal rational set $X_{1} \subseteq S$. By
Lemma~\ref{nozero}, we may assume that $0 \notin X_0$ and $0 \notin
X_1$.

Let $C$ and $D$ be automata over $S$ accepting $X_0$ and $X_1$ respectively.
Since $C$, $D$ and $A$ have only finitely many edges between them,
we may choose finite subsets $I' \subseteq I$ and $J' \subseteq J$
such that the edge labels of $C$ and $D$ all lie in $I' \times G
\times J'$, and the edge labels of $A$ all lie in $(I' \times G
\times J') \times \Sigma^*$.

For each $i \in I'$ and $j \in J'$, we let $X_{ij} = \{g \in G \mid
(i,g,j) \in X_0 \}$.  By Proposition \ref{subset}, each $X_{ij}$ is
a rational subset of $G$. It follows that
$$X_{ij}' \ = \ X_{ij} \times \lbrace \epsilon \rbrace$$
is a rational subset of $G \times \Sigma^*$; let $C_{ij}$ be an
automaton accepting $X_{ij}'$.

Similarly, for each $i \in I'$ and $j \in J'$ we define $Y_{ij} =
\{g^{-1} \in G \mid (i,g,j) \in X_{1} \}$. By
Propositions~\ref{subset} and~\ref{ratinv}, $Y_{ij}$ is a rational
subset of $G$, and so
$$Y_{ij}' \ = \ Y_{ij} \times \lbrace \epsilon \rbrace$$
is a rational subset of $G \times \Sigma^*$; let $D_{ij}$ be an
automaton accepting $Y_{ij}'$.

Assume without loss of generality that the automaton $A$ and all the
automata $C_{ij}$ and $D_{ij}$ have disjoint vertex sets. We
construct from these automata a $G$-automaton $B$ with
\begin{itemize}
\item vertex set the union of the vertex sets of $C_{ij}$ and $D_{ij}$
(for $i \in I'$ and $j \in J'$) together with $I' \times Q \times
J'$ where $Q$ is the vertex set of $A$, and a new vertex $q_0'$;
\item initial vertex $q_0'$;
\item terminal vertices the terminal vertices of the automata $D_{ij}$;
\item all the edges of the automata $C_{ij}$ and $D_{ij}$;
\item for each $i \in I'$ and $j \in J'$, an edge from $q_0'$ to the
initial vertex of $C_{ij}$ labelled $(1, \epsilon)$;
\item for each $i \in I'$ and $j \in J'$, an edge from each terminal vertex
of $C_{ij}$ to $(i, q_0, j)$ labelled $(1, \epsilon)$, where $q_0$ is
the initial vertex for $A$;
\item for each edge in $A$ from $p$ to $q$ labelled $((i, g, j), w)$ and
each $i' \in I'$ and $j' \in J'$, an edge from $(i', p, j')$ to
$(i', q, j)$ labelled $(P_{j' i} g, w)$;
\item for each $i \in I'$, $j \in J'$ and terminal vertex $p$ of $A$, an
edge from $(i, p, j)$ to the initial vertex of $D_{ij}$ labelled
 $(1, \epsilon)$.
\end{itemize}
Since $I'$, $J'$ and all the automata $A$, $C_{ij}$ and $D_{ij}$ are
finite, it follows that the $G$-automaton $B$ is finite. We now show
that $B$ accepts the language $L$.

Let \mbox{$w \in L$}. Then there exists a path through the automaton
$A$ labelled $((i,g,j),w)$ connecting the initial vertex with some
terminal vertex ($p_{t}$ say), such that \[(i_{0},g_{0},j_{0})(i,g,j)
= (i',g',j') \in X_{1}\] for some $(i_{0},g_{0},j_{0}) \in X_{0}$.
Suppose this path has the form
\[
\mathbf{q_{0}}\xrightarrow{((i_1,g_1,j_1),
x_1)}\mathbf{q_{1}}\xrightarrow{((i_{2},g_{2},j_{2}),
x_2)}\mathbf{q_{2}}\xrightarrow{((i_{3},g_{3},j_{3}), x_3)} \ldots
\mathbf{q_{m-1}}\xrightarrow{((i_{m},g_{m},j_m), x_m)}\mathbf{q_{m}}
\]
where $q_0$ is the initial vertex and $q_m = p_{t}$ is a terminal
vertex of $A$ and $w = x_1 \dots x_m$. Note that we must have $i' = i_0$,
$j' = j_m$ and
$$g = g_1 P_{j_1 i_2} g_2 \dots P_{j_{m-1} i_m} g_m.$$
Now by construction, $B$ has a path $\pi_2$ of the form
\begin{align*}
\mathbf{(i_0, q_{0}, j_0)}\xrightarrow{(P_{j_0 i_1} g_1,
x_1)}\mathbf{(i_0, q_1, j_1)}\xrightarrow{(P_{j_1 i_2} g_2,
x_2)}\mathbf{(i_0,q_{2},j_2)}\xrightarrow{(P_{j_2 i_3} g_3, x_3)}
\ldots \\
\ldots \xrightarrow{(P_{j_{m-1} i_m}g_{m},
x_m)}\mathbf{(i_0,q_m,j_m)}
\end{align*}
Moreover, from the fact that $(i_0, g_0, j_0) \in X_0$ we see that
$g_0 \in X_{i_0j_0}$, so that $(g_0,\epsilon) \in X_{i_0 j_0}'$.
Hence, $(g_0, \epsilon)$ labels a path in $C_{i_0j_0}$ from the initial
vertex to a terminal vertex. It follows easily that $(g_0,
\epsilon)$ labels a path $\pi_1$ in $B$ from the initial vertex
$q_0'$ to $(i_0, q_0, j_0)$ where $q_0$. Similarly, since $(i', g',
j') \in X_1$ we deduce that $((g')^{-1}, \epsilon) \in Y_{i' j'} = Y_{i_0 j_m}$
so that $B$ has a path $\pi_3$ from $(i_0, q_m, j_m)$ to a terminal
vertex labelled $((g')^{-1}, \epsilon)$.

Composing the paths $\pi_1$, $\pi_2$ and $\pi_3$, we see that $B$
has a path from the initial vertex to a terminal vertex with label
$$(g_0 P_{j_0 i_1} g_1 P_{j_1 i_2} g_2 \dots P_{j_{m-1} i_m} g_m (g')^{-1}, x_1 x_2 \dots x_m)$$
But we know that $(i_{0},g_{0},j_{0})(i,g,j) = (i',g',j')$, so we
must have
$$g_0 P_{j_0 i_1} g_1 P_{j_1 i_2} g_2 \dots P_{j_{m-1} i_m} g_m = g'$$
and hence
$$g_0 P_{j_0 i_1} g_1 P_{j_1 i_2} g_2 \dots P_{j_{m-1} i_m} g_m (g')^{-1} = 1.$$
It follows that $w$ is accepted by the $G$-automaton $B$, as required.

Conversely, suppose $w$ is accepted by the $G$-automaton $B$. Then
there is a path in $B$ from the initial vertex to a terminal vertex
labelled $(1,w)$. We deduce easily from the construction of $B$ that this
path must have the form $\pi_1 \pi_2 \pi_3$ where
\begin{itemize}
\item $\pi_1$ runs from the start vertex to some vertex $(i_0, q_0, j_0)$ with
label of the form $(g_0, \epsilon)$ for some $g_0 \in X_{i_0 j_0}$,
so that $(i_0, g_0, j_0) \in X_0$;
\item $\pi_2$ runs from $(i_0, q_0, j_0)$ to a vertex $(i_0, q_m, j_m)$
where $q_m$ is a terminal vertex of $A$; and
\item $\pi_3$ runs from $(i_0, q_m, j_m)$ to a terminal vertex with label
$((g')^{-1}, \epsilon)$ where $(g')^{-1} \in Y_{i_0 j_m}$, so that
$(i_0, g', j_m) \in X_1$.
\end{itemize}
Moreover, $\pi_2$ must have the form
\begin{align*}
\mathbf{(i_0, q_{0}, j_0)}\xrightarrow{(P_{j_0 i_1} g_1,
x_1)}\mathbf{(i_0, q_1, j_1)}\xrightarrow{(P_{j_1 i_2} g_2,
x_2)}\mathbf{(i_0,q_{2},j_2)}\xrightarrow{(P_{j_2 i_3} g_3, x_3)}
\ldots \\
\ldots \xrightarrow{(P_{j_{m-1},i_m}g_{m},
x_m)}\mathbf{(i_0,q_m,j_m)}
\end{align*}
where, since the label of the entire path $\pi$ is $(1, w)$, we must
have $w = x_1 \dots x_m$ and $g_0 P_{j_0 i_1} g_1 \dots P_{j_{m-1}
i_m} g_m (g')^{-1} = 1$, that is,
$$g_0 P_{j_0 i_1} g_1 \dots P_{j_{m-1} i_m} g_m = g'.$$
We deduce from the path above and the construction of $B$ that $A$ has a path
\[
\mathbf{q_{0}}\xrightarrow{((i_1,g_1,j_1),
x_1)}\mathbf{q_{1}}\xrightarrow{((i_{2},g_{2},j_{2}),
x_2)}\mathbf{q_{2}}\xrightarrow{((i_{3},g_{3},j_{3}), x_3)} \ldots
\mathbf{q_{m-1}}\xrightarrow{((i_{m},g_{m},j_m), x_m)}\mathbf{q_{m}}
\]
Since $q_0$ and $q_m$ are initial and terminal vertices of $A$
respectively, it follows that $A$ accepts $(x, w)$ where
$$x = (i_1, g_1, j_1) (i_2, g_2, j_2) \dots (i_m, g_m, j_m).$$
But $(i_0, g_0, j_0)$ lies in $X_0$ and
\begin{align*}
(i_0, g_0, j_0) x &= (i_0, g_0, j_0) (i_1, g_1, j_1) \dots (i_m, g_m, j_m) \\
&= (i_0, g_0 P_{j_0 i_1} g_1 \dots P_{j_{m-1} i_m} g_m, j_m) \\
&= (i_0, g', j_m)
\end{align*}
lies in $X_1$, from which we deduce that the rational $S$-automaton
$A$ accepts the word $w$, and so $w \in L$ as required.
\end{proof}

\section*{Acknowledgements}

The research of the second author was supported by an RCUK Academic
Fellowship.

\bibliographystyle{plain}

\end{document}